# The Geometric Proof of the Hecke Conjecture

**Kaida Shi**

(Department of Mathematics, Zhejiang Ocean University,

Zhoushan City, 316004, Zhejiang Province, China)

**Abstract** Beginning from the resolution of Dirichlet *L* function $L(s,\ \chi)$, using the inner product formula of two infinite-dimensional vectors in the complex space, the author proved the baffling problem---- Hecke conjecture.

**Keywords: Dirichlet *L* function $L(s,\ \chi)$, Hecke conjecture, non-trivial zeroes, solid of rotation, axis-cross section, bary-centric coordinate, inner product, infinite-dimensional vectors.**

## 1 Introduction

For investigating the distribution problem of prime numbers within the arithmetic series, first, L. Dirichlet imported the Dirichlet *L* function $L(s,\ \chi)$. Although its property and effect is similar with Riemann Zeta function $\zeta(s)$, but the different matter is, the research about the distribution of the real zeros of the Dirichlet *L* function $L(s,\ \chi)$ which correspond real property is a very difficult thing. But because of this, it has very important meaning.

For solving above-mentioned problem, Hecke put forward a famous proposition as follows:

**Hecke Conjecture:** Suppose that $\chi$ is a real property, then we have

$$L(s,\ \chi) \neq 0,\ \ 0 < s < 1.$$

Now, let's prove the proposition.

## 2 The resolution of Dirichlet *L* function $L(s,\chi)$

Considering the geometric meaning of the Dirichlet *L* function

$$L(s,\chi) = \sum_{n=1}^{\infty} \frac{\chi(n)}{n^s} = \frac{\chi(1)}{1^s} + \frac{\chi(2)}{2^s} + \frac{\chi(3)}{3^s} + \cdots + \frac{\chi(n)}{n^s} + \cdots, \tag{1}$$

we have the following figure:

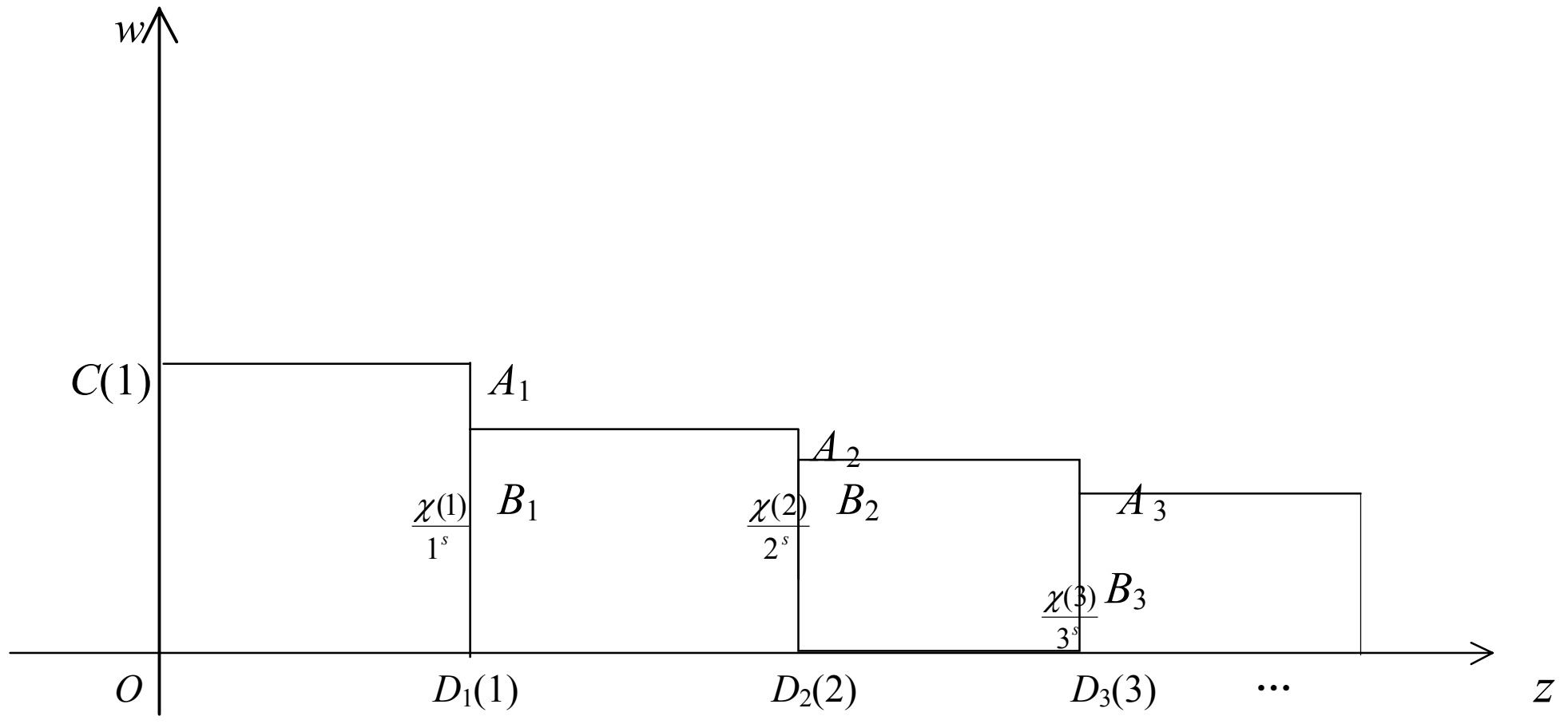

Fig. 2

In this figure, the areas of the rectangles $A_1D_1OC, A_2D_2D_1B_1, A_3D_3D_2B_2, \cdots$ are respectively:

$$\frac{\chi(1)}{1^s}\cdot 1,\ \frac{\chi(2)}{2^s}\cdot 1,\ \frac{\chi(3)}{3^s}\cdot 1,\ \cdots,\ \frac{\chi(n)}{n^s}\cdot 1,\ \cdots$$

therefore, the geometric meaning of the Dirichlet $L$ function $L(s,\chi)$ is the sum of the areas of a series of rectangles within the complex space $s$.

Using the inner product formula between two infinite-dimensional vectors, the Dirichlet $L$ function $L(s,\chi)$ equation

$$L(s,\chi)=\frac{\chi(1)}{1^s}+\frac{\chi(2)}{2^s}+\frac{\chi(3)}{3^s}+\cdots+\frac{\chi(n)}{n^s}+\cdots=0, \qquad (s=\sigma+it)$$

can be resolved as

$$L(s,\chi)=(\frac{\chi(1)}{1^\sigma},\frac{\chi(2)}{2^\sigma},\frac{\chi(3)}{3^\sigma},\cdots,\frac{\chi(n)}{n^\sigma},\cdots)\cdot(\frac{1}{1^{it}},\frac{1}{2^{it}},\frac{1}{3^{it}},\cdots,\frac{1}{n^{it}},\cdots)=0. \qquad (2)$$

or

$$L(s,\chi)=(\frac{1}{1^\sigma},\frac{1}{2^\sigma},\frac{1}{3^\sigma},\cdots,\frac{1}{n^\sigma},\cdots)\cdot(\frac{\chi(1)}{1^{it}},\frac{\chi(2)}{2^{it}},\frac{\chi(3)}{3^{it}},\cdots,\frac{\chi(n)}{n^{it}},\cdots)=0\,. \qquad (3)$$

From the expressions (2) and (3), we obtain:

$$\sqrt{\sum_{n=1}^{\infty}\frac{\chi^2(n)}{n^{2\sigma}}}\cdot\sqrt{\sum_{n=1}^{\infty}\frac{1}{n^{2it}}}\cdot\cos(\widehat{\vec{N}_1,\vec{N}_2})=0\,; \tag{4}$$

and

$$\sqrt{\sum_{n=1}^{\infty}\frac{1}{n^{2\sigma}}}\cdot\sqrt{\sum_{n=1}^{\infty}\frac{\chi^2(n)}{n^{2it}}}\cdot\cos(\widehat{\vec{N}_3,\vec{N}_4})=0\,. \tag{5}$$

But in the complex space, if the inner product between two vectors equals to zero, then these two vectors are perpendicular, namely, $(\widehat{\vec{N}_1,\vec{N}_2})=\frac{\pi}{2}$, $(\widehat{\vec{N}_3,\vec{N}_4})=\frac{\pi}{2}$, therefore, we have

$$\cos(\widehat{\vec{N}_1,\vec{N}_2})=0,$$

and

$$\cos(\widehat{\vec{N}_3,\vec{N}_4})=0.$$

## 3 The relationship between the volume of the rotation solid and the area of its axis-cross section within the complex space

From the expression (1), we can know that when $\mathrm{Re}(s)=\frac{1}{2}$, the series (harmonic series) within the first radical expression of left side is divergent. Because $\cos(\widehat{\vec{N}_1,\vec{N}_2})=0,$ therefore the situation of the complex series within the second radical expression of left side will change unable to research. Hence, we must transform the Riemann Zeta function $\zeta(s)$ equation. For this aim, let's derive the relationship between the volume of the rotation solid and the area of its axis-cross section within the complex space.

We call the cross section which pass through the axis $z$ and intersects the rotation solid as **axis-cross section**. According to the barycentric formula of the complex plane lamina:

$$\begin{cases} \xi = \dfrac{\int_a^b z[f(z)-g(z)]dz}{\int_a^b [f(z)-g(z)]dz}, \\ \eta = \dfrac{\frac{1}{2}\int_a^b [f^2(z)-g^2(z)]dz}{\int_a^b [f(z)-g(z)]dz}. \end{cases}$$

we have

$$2\pi\eta = \frac{\pi\int_a^b [f^2(z)-g^2(z)]dz}{\int_a^b [f(z)-g(z)]dz}. \qquad (6)$$

The numerator of the fraction of right side of the formula (**) is the volume of the rotation solid, and the denominator is the area of the axis-cross section. The $\eta$ of left side of the formula (6) is the longitudinal coordinate of the barycenter of axis-cross section. Its geometric explanation is: taking the axis-cross section around the axis $z$ to rotate the angle of $2\pi$, we will obtain the volume of the rotation solid.

## 4 The proof of the Hecke conjecture

Because we have the volume formula of the rotation solid formed by rotating the rectangle $A_n D_n D_{n-1} B_{n-1}$ around the axis $z$:

$$V_n = \int_{n-1}^{n} \frac{\pi\chi^2(n)}{n^{2s}} dz = \frac{\pi\chi^2(n)}{n^{2s}} z\Big|_{n-1}^{n} = \frac{\pi\chi^2(n)}{n^{2s}},$$

therefore, the sum of the volumes of a series of the cylinders formed by rotating a series of the rectangles $A_1D_1OC, A_2D_2D_1B_1, A_3D_3D_2B_2, \cdots$ around the axis $z$ is:

$$V = \sum_{n=1}^{\infty} V_n = \pi\left(\frac{\chi^2(1)}{1^{2s}} + \frac{\chi^2(2)}{2^{2s}} + \frac{\chi^2(3)}{3^{2s}} + \cdots + \frac{\chi^2(n)}{n^{2s}} + \cdots\right).$$

But the sum of the areas of a series of the rectangles $A_1D_1OC, A_2D_2D_1B_1, A_3D_3D_2B_2, \cdots$ is:

$$\frac{\chi(1)}{1^s}+\frac{\chi(2)}{2^s}+\frac{\chi(3)}{3^s}+\cdots+\frac{\chi(n)}{n^s}+\cdots.$$

By the formula (6), we have

$$\pi\left(\frac{\chi^2(1)}{1^{2s}}+\frac{\chi^2(2)}{2^{2s}}+\frac{\chi^2(3)}{3^{2s}}+\cdots+\frac{\chi^2(n)}{n^{2s}}+\cdots\right)$$

$$=2\pi\eta\left(\frac{\chi(1)}{1^s}+\frac{\chi(2)}{2^s}+\frac{\chi(3)}{3^s}+\cdots+\frac{\chi(n)}{n^s}+\cdots\right). \tag{7}$$

Substituting the Dirichlet $L$ function $L(s,\chi)$ equation

$$L(s,\chi)=\frac{\chi(1)}{1^s}+\frac{\chi(2)}{2^s}+\frac{\chi(3)}{3^s}+\cdots+\frac{\chi(n)}{n^s}+\cdots=0$$

into (7), we can obtain the transformed equation:

$$\frac{\chi^2(1)}{1^{2s}}+\frac{\chi^2(2)}{2^{2s}}+\frac{\chi^2(3)}{3^{2s}}+\cdots+\frac{\chi^2(n)}{n^{2s}}+\cdots=0.$$

This can be explained geometrically as: **if the area of the axis-cross section equals to zero, then it corresponds to the volume of the rotation solid also equals to zero**.

Because $s=\sigma+0i$, therefore, above expression becomes

$$\frac{\chi^2(1)}{1^{2\sigma}}+\frac{\chi^2(2)}{2^{2\sigma}}+\frac{\chi^2(3)}{3^{2\sigma}}+\cdots+\frac{\chi^2(n)}{n^{2\sigma}}+\cdots=0.$$

According to the inner product formula between two infinite-dimensional vectors, above expression can be written as:

$$\sqrt{\frac{\chi^4(1)}{1^{4\sigma}}+\frac{\chi^4(2)}{2^{4\sigma}}+\frac{\chi^4(3)}{3^{4\sigma}}+\cdots+\frac{\chi^4(n)}{n^{4\sigma}}+\cdots}\cdot\sqrt{\frac{1}{1^{4ti}}+\frac{1}{2^{4ti}}+\frac{1}{3^{4ti}}+\cdots+\frac{1}{n^{4ti}}+\cdots}\cdot\cos(\widehat{\overrightarrow{N_1},\overrightarrow{N_2}})=0,$$

In the expression, when $t=0$, although the series

$$\sum_{n=1}^{\infty}\sin(0\cdot\ln n)=0+0+0+\cdots+0+\cdots=0,$$

but because the series

$$\sum_{n=1}^{\infty}\cos(0\cdot\ln n)=1+1+1+\cdots+1+\cdots=\infty,$$

therefore, the complex series

$$\begin{aligned}&\frac{1}{1^{0i}}+\frac{1}{2^{0i}}+\frac{1}{3^{0i}}+\cdots+\frac{1}{n^{0i}}+\cdots\\&=\sum_{n=1}^{\infty}[\cos(0\ln n)-i\sin(0\ln n)]\\&=(1+1+1+\cdots+1+\cdots)-\\&-i\cdot(0+0+0+\cdots+0+\cdots)=\infty.\end{aligned}$$

Because "infinite" multiplied by zero doesn't equal to zero, therefore

$$\begin{aligned}&\sqrt{\sum_{n=1}^{\infty}\frac{\chi^4(n)}{n^{4\sigma}}}\cdot\sqrt{\sum_{n=1}^{\infty}\frac{1}{n^{0i}}}\cdot 0\\&=\sqrt{\sum_{n=1}^{\infty}\frac{\chi^4(n)}{n^{4\sigma}}}\cdot\sqrt{\infty}\cdot 0\neq 0.\end{aligned}$$

On the other hand, from the expression (3), we have:

$$\sqrt{\frac{1}{1^{4\sigma}}+\frac{1}{2^{4\sigma}}+\frac{1}{3^{4\sigma}}+\cdots+\frac{1}{n^{4\sigma}}+\cdots}\cdot\sqrt{\frac{\chi^4(1)}{1^{4it}}+\frac{\chi^4(2)}{2^{4it}}+\frac{\chi^4(3)}{3^{4it}}+\cdots+\frac{\chi^4(n)}{n^{4it}}+\cdots}\cdot\cos(\widehat{\overrightarrow{N_3},\overrightarrow{N_4}})=0,$$

As is known to all, among the series:

$$l,\ l+q,\ l+2q,\ \cdots\cdots\quad (l,\ q)=1.$$

we can know $\chi(n;\ q)=1$. According to the character of the property, we have

$$\chi(n;\ q)=\chi(n+q;\ q),$$

also, we have

$$\chi(n+q;\ q)=\chi(n+2q;\ q)=\chi(n+3q;\ q)=\cdots\cdots=1,$$

therefore, the real properties are infinite. Hence, we can obtain:

$$\sum_{n=1}^{\infty}\frac{\chi^4(n)}{n^{0i}}=\sum_{n=1}^{\infty}\chi^4(n)(\cos(0\cdot\log n)-i\sin(0\cdot\log n))$$
$$=1+\cdots+1+\cdots+1+\cdots+1+\cdots=\infty.$$

therefore, we have

$$\sqrt{\sum_{n=1}^{\infty}\frac{1}{n^{4\sigma}}}\cdot\sqrt{\sum_{n=1}^{\infty}\frac{\chi^4(n)}{n^{0i}}}\cdot 0$$
$$=\sqrt{\sum_{n=1}^{\infty}\frac{1}{n^{4\sigma}}}\cdot\sqrt{\infty}\cdot 0\neq 0.$$

So, we have proved that when the $\chi$ is real property, then we have

$$L(s,\ \chi)\neq 0,\ 0<s<1.$$

# Appendix

## The inner product formula between two vectors in real space can be extended formally to complex space

**Dear Mr. Referee,**

**Thank you for your review, please observe following examples.**

**Example 1** $A\{1+i,\ 3\},\ B\{-i,\ 2i\},\ C\{1,\ -i\}$.

$\overrightarrow{AB}=\{-1-2i,\ -3+2i\},\ \overrightarrow{AC}=\{-i,\ -3-i\},\ \overrightarrow{BC}=\{1+i,\ -3i\}$.

$\left|\overrightarrow{AB}\right|=\sqrt{(-1-2i)^2+(-3+2i)^2}=\sqrt{2-8i}$,

$\left|\overrightarrow{AC}\right|=\sqrt{(-i)^2+(-3-i)^2}=\sqrt{7+6i}$,

$\left|\overrightarrow{BC}\right|=\sqrt{(1+i)^2+(-3i)^2}=\sqrt{-9+2i}$.

According to **Cosine theorem**

$$\left|\overrightarrow{AB}\right|^2+\left|\overrightarrow{AC}\right|^2-\left|\overrightarrow{BC}\right|^2=2\left|\overrightarrow{AB}\right|\left|\overrightarrow{AC}\right|\cos(\overrightarrow{AB},\overrightarrow{AC}),$$

we have

$$\frac{\left|\overrightarrow{AB}\right|^2+\left|\overrightarrow{AC}\right|^2-\left|\overrightarrow{BC}\right|^2}{2}=\frac{(2-8i)+(7+6i)-(-9+2i)}{2}=9-2i.$$

On the other hand, we have

$$\overrightarrow{AB}\cdot\overrightarrow{AC}=\left|\overrightarrow{AB}\right|\left|\overrightarrow{AC}\right|\cos(\overrightarrow{AB},\overrightarrow{AC})=(-1-2i)(-i)+(-3+2i)(-3-i)=9-2i.$$

According to the **Cosine theorem**

$$\left|\overrightarrow{AC}\right|^2+\left|\overrightarrow{BC}\right|^2-\left|\overrightarrow{AB}\right|^2=2\left|\overrightarrow{AC}\right|\left|\overrightarrow{BC}\right|\cos(\overrightarrow{AC},\overrightarrow{BC})$$

$$\frac{\left|\overrightarrow{AC}\right|^2+\left|\overrightarrow{BC}\right|^2-\left|\overrightarrow{AB}\right|^2}{2}=\frac{(7+6i)+(-9+2i)-(2-8i)}{2}=-2+8i.$$

On the other hand, we have

$$\overrightarrow{AC}\cdot\overrightarrow{BC}=\left|\overrightarrow{AB}\right|\left|\overrightarrow{AC}\right|\cos(\overrightarrow{AB},\overrightarrow{AC})=(-i)(1+i)+(-3-i)(-3i)=-2+8i.$$

$$S_{\Delta ABC}=\frac{1}{2}\left|\overrightarrow{AB}\right|\left|\overrightarrow{AC}\right|\sin\theta_1=\frac{1}{2}\sqrt{2-8i}\cdot\sqrt{7+6i}\cdot\frac{\sqrt{-15-8i}}{\sqrt{62-44i}}=\frac{1}{2}\sqrt{-15-8i},$$

$$S_{\Delta ACB}=\frac{1}{2}\left|\overrightarrow{AC}\right|\left|\overrightarrow{BC}\right|\sin\theta_{21}=\frac{1}{2}\sqrt{7+6i}\cdot\sqrt{-9+2i}\cdot\frac{\sqrt{-15-8i}}{\sqrt{-75-40i}}=\frac{1}{2}\sqrt{-15-8i}.$$

**Example 2** $A\{1+i,\ 1-i,\ 2i\},\ B\{1-i,\ 1+i,\ -2i\},\ C\{1,\ 0,\ i\}.$

$\overrightarrow{AB}=\{-2i,\ 2i,\ -4i\},\ \overrightarrow{AC}=\{-i,\ -1+i,\ -i\},\ \overrightarrow{BC}=\{i,\ -1-i,\ 3i\}.$

$\left|\overrightarrow{AB}\right|=\sqrt{(-2i)^2+(2i)^2+(-4i)^2}=\sqrt{-24},$

$\left|\overrightarrow{AC}\right|=\sqrt{(-i)^2+(-1+i)^2+(-i)^2}=\sqrt{-2-2i},$

$\left|\overrightarrow{BC}\right|=\sqrt{(i)^2+(-1-i)^2+(3i)^2}=\sqrt{-10+2i}.$

According to the **Cosine theorem**

$$\left|\overrightarrow{AB}\right|^2+\left|\overrightarrow{AC}\right|^2-\left|\overrightarrow{BC}\right|^2=2\left|\overrightarrow{AB}\right|\left|\overrightarrow{AC}\right|\cos(\overrightarrow{AB},\overrightarrow{AC}),$$

we have

$$\frac{\left|\overrightarrow{AB}\right|^2+\left|\overrightarrow{AC}\right|^2-\left|\overrightarrow{BC}\right|^2}{2}=\frac{(-24)+(-2-2i)-(-10+2i)}{2}=-8-2i.$$

On the other hand, we have

$$\overrightarrow{AB}\cdot\overrightarrow{AC}=\left|\overrightarrow{AB}\right|\left|\overrightarrow{AC}\right|\cos(\overrightarrow{AB},\overrightarrow{AC})=(-2i)(-i)+(2i)(-1+i)+(-4i)(-i)=-8-2i.$$

According to the **Cosine theorem**

$$\left|\overrightarrow{AC}\right|^2+\left|\overrightarrow{BC}\right|^2-\left|\overrightarrow{AC}\right|^2=2\left|\overrightarrow{AC}\right|\left|\overrightarrow{BC}\right|\cos(\overrightarrow{AC},\overrightarrow{BC}),$$

we have

$$\frac{\left|\overrightarrow{AC}\right|^2+\left|\overrightarrow{BC}\right|^2-\left|\overrightarrow{AB}\right|^2}{2}=\frac{(-2-2i)+(-10+2i)-(-24)}{2}=6.$$

On the other hand, we have

$$\overrightarrow{AC}\cdot\overrightarrow{BC}=\left|\overrightarrow{AC}\right|\left|\overrightarrow{BC}\right|\cos(\overrightarrow{AC},\overrightarrow{BC})=(-i)(i)+(-1+i)(-1-i)+(-i)(3i)=6.$$

$$S_{\Delta ABC}=\frac{1}{2}\left|\overrightarrow{AB}\right|\left|\overrightarrow{AC}\right|\sin\theta_1=\frac{1}{2}\sqrt{-24}\cdot\sqrt{-2-2i}\cdot\frac{\sqrt{-12+16i}}{\sqrt{48+48i}}=\sqrt{-3+4i},$$

$$S_{\Delta ACB}=\frac{1}{2}\left|\overrightarrow{AC}\right|\left|\overrightarrow{BC}\right|\sin\theta_{21}=\frac{1}{2}\sqrt{-2-2i}\cdot\sqrt{-10+2i}\cdot\frac{\sqrt{-12+16i}}{\sqrt{24+16i}}=\sqrt{-3+4i}.$$

**Example 3** $A\{8i,\ 14,\ 8-i,\ 1\}$, $B\{6,\ 15i,\ 17,\ -8\}$, $C\{3-i,\ 10+7i,\ 11,\ 3i\}$.

$\overrightarrow{AB}=\{6-8i,\ -14+15i,\ 9+i,\ -9\}$, $\overrightarrow{AC}=\{3-9i,\ -4+7i,\ 3+i,\ -1+3i\}$,

$\overrightarrow{BC}=\{-3-i,\ 10-8i,\ -6,\ 8+3i\}$.

$\left|\overrightarrow{AB}\right|=\sqrt{(6-8i)^2+(-14+5i)^2+(9+i)^2+(-9)^2}=\sqrt{104-498i}$,

$\left|\overrightarrow{AC}\right|=\sqrt{(3-9i)^2+(-4+7i)^2+(3+i)^2+(-1+3i)^2}=\sqrt{-105-110i}$,

$\left|\overrightarrow{BC}\right|=\sqrt{(-3-i)^2+(10-8i)^2+(-6)^2+(8+3i)^2}=\sqrt{135-106i}$.

According to the **Cosine theorem**

$$\left|\overrightarrow{AB}\right|^2+\left|\overrightarrow{AC}\right|^2-\left|\overrightarrow{BC}\right|^2=2\left|\overrightarrow{AB}\right|\left|\overrightarrow{AC}\right|\cos(\overrightarrow{AB},\overrightarrow{AC}),$$

we have

$$\frac{\left|\overrightarrow{AB}\right|^2+\left|\overrightarrow{AC}\right|^2-\left|\overrightarrow{BC}\right|^2}{2}=\frac{(104-498i)+(-105-110i)-(-135-106i)}{2}=-68-251i.$$

On the other hand, we have

$$\overrightarrow{AB}\cdot\overrightarrow{AC}=\left|\overrightarrow{AB}\right|\left|\overrightarrow{AC}\right|\cos(\overrightarrow{AB},\overrightarrow{AC})=$$
$$=(6-8i)(3-9i)+(-14+15i)(-4+7i)+(9+i)(3+i)+(-9)(-1+3i)=-68-251i.$$

According to the **Cosine theorem**

$$\left|\overrightarrow{AC}\right|^2+\left|\overrightarrow{BC}\right|^2-\left|\overrightarrow{AC}\right|^2=2\left|\overrightarrow{AC}\right|\left|\overrightarrow{BC}\right|\cos(\overrightarrow{AC},\overrightarrow{BC}),$$

we have

$$\frac{\left|\overrightarrow{AC}\right|^2+\left|\overrightarrow{BC}\right|^2-\left|\overrightarrow{AB}\right|^2}{2}=\frac{(-105-110i)+(135-106i)-(104-498i)}{2}=-37+141i.$$

On the other hand, we hand

$$\overrightarrow{AC}\cdot\overrightarrow{BC}=\left|\overrightarrow{AC}\right|\left|\overrightarrow{BC}\right|\cos(\overrightarrow{AC},\overrightarrow{BC})=$$
$$=(3-9i)(-3-i)+(-4+7i)(10-8i)+(-6)(3+i)+(-1+3i)(8+3i)=-37+141i.$$

$$S_{\Delta ABC}=\frac{1}{2}\left|\overrightarrow{AB}\right|\left|\overrightarrow{AC}\right|\sin\theta_1=\frac{1}{2}\sqrt{104-498i}\cdot\sqrt{-105-110i}\cdot\frac{\sqrt{-7323+6714i}}{\sqrt{-65700+40850i}}$$
$$=\frac{1}{2}\sqrt{-7323+6714i},$$

$$S_{\Delta ACB}=\frac{1}{2}\left|\overrightarrow{AC}\right|\left|\overrightarrow{BC}\right|\sin\theta_{21}=\frac{1}{2}\sqrt{-105-110i}\cdot\sqrt{135-106i}\cdot\frac{\sqrt{-7323+6714i}}{\sqrt{-25835-3720i}}$$
$$=\frac{1}{2}\sqrt{-7323+6714i}.$$